\documentclass{amsart}[11]

\usepackage{amsfonts}
\usepackage{amssymb}
\usepackage{amsthm}

\usepackage{amsmath}
\usepackage{amscd}
\usepackage{enumerate}
\usepackage{graphicx}
\usepackage[all]{xy}
\usepackage{mathrsfs}

\theoremstyle{plain}
\newtheorem{thm}{Theorem}
\newtheorem{prop}[thm]{Proposition}
\newtheorem{lem}[thm]{Lemma}

\newtheorem{defn}[thm]{Definition}
\newtheorem{cor}[thm]{Corollary}

\theoremstyle{remark}
\newtheorem{rem}[thm]{Remark}
\newtheorem{eg}[thm]{Example}

\newcommand{\bra}{{\langle}}
\newcommand{\ket}{{\rangle}}
\newcommand{\dd}{{\partial}}

\newcommand{\bA}{\mathbb{A}}
\newcommand{\bC}{\mathbb{C}}

\newcommand{\bP}{\mathbb{P}}

\newcommand{\bZ}{\mathbb{Z}}

\newcommand{\cM}{{\mathcal M}}

\newcommand{\cO}{{\mathcal O}}

\newcommand{\cU}{{\mathcal U}}

\newcommand{\sC}{\mathscr{C}}
\newcommand{\sD}{\mathscr{D}}
\newcommand{\sE}{\mathscr{E}}
\newcommand{\sF}{\mathscr{F}}
\newcommand{\sG}{\mathscr{G}}

\newcommand{\sL}{\mathscr{L}}

\newcommand{\sO}{\mathscr{O}}

\newcommand{\sR}{\mathscr{R}}

\newcommand{\sV}{\mathscr{V}}

\newcommand{\fh}{\mathfrak{h}}

\newcommand{\fsl}{\mathfrak{sl}}
\newcommand{\fg}{\mathfrak{g}}
\newcommand{\fl}{\mathfrak{l}}
\newcommand{\fn}{\mathfrak{n}}
\newcommand{\fp}{\mathfrak{p}}

\newcommand{\fu}{\mathfrak{u}}

\DeclareMathOperator{\Ad}  {Ad}
\DeclareMathOperator{\Aut}  {Aut}
\DeclareMathOperator{\cok} {cok}
\DeclareMathOperator{\dbar}{\bar{\partial}}

\DeclareMathOperator{\End} {End}
\DeclareMathOperator{\Ext} {Ext}
\DeclareMathOperator{\GL}  {GL}

\DeclareMathOperator{\Gr}{Gr}

\DeclareMathOperator{\Hom} {Hom}

\DeclareMathOperator{\im}  {im}
\DeclareMathOperator{\Jac} {Jac}

\DeclareMathOperator{\Pic} {Pic}

\DeclareMathOperator{\SL}  {SL}
\DeclareMathOperator{\Spec}{Spec}
\DeclareMathOperator{\Sym} {Sym}
\DeclareMathOperator{\Tr}  {Tr}

\begin{document}

\title{Poisson geometry of parabolic bundles on elliptic curves}
\author{David Balduzzi}
\address{
Department of Psychiatry, University of Wisconsin, Madison, WI 53705 \\
{\em E-mail:} {\tt balduzzi@wisc.edu}
}

\maketitle

\begin{abstract}
The moduli space of $G$-bundles on an elliptic curve with additional flag structure admits a Poisson structure. The bivector can be defined using double loop group, loop group and sheaf cohomology constructions. We investigate the links between these methods and for the case $\SL_2$ perform explicit computations, describing the bracket and its leaves in detail.
\end{abstract}


\section{Introduction}

Let $E$ be an elliptic curve. In this note we define a natural Poisson structure on the moduli space of parabolic bundles on $E$. The construction is inspired by a Poisson structure constructed by Polishchuk \cite{polishchuk:98} on the moduli space of stable triples $(E_1,E_2,\Phi)$ where $\Phi:E_1\rightarrow E_2$ is a morphism between vector bundles over $E$. Less directly the idea comes from Mukai's \cite{muk:84} construction of a Poisson structure on moduli space of sheaves on abelian surfaces. In our case the Poisson bracket can be further motivated as a Hamiltonian reduction of a Kirillov-Kostant bracket on the dual to the Lie algebra of a double loop group.

We study the bracket in detail in the case of $\SL_2$ in particular examining its leaves. We also provide a third construction of Poisson brackets on the moduli space using loop groups and $r$-matrices, and investigate the relationship between the two brackets for $\SL_2$, where the computations are tractable.

The contents of the note are as follows. In the next section we collect some facts due to Friedman-Morgan-Witten on moduli spaces of principal bundles on elliptic curves. Section three constructs the Poisson bivector using sheaf cohomology, and finds some Casimirs using group-theoretic maps. It also provides a second construction of the Poisson bracket using Atiyah-Bott reduction and double loop groups, motivating the first. The idea is that the affine space of complex structures on a principal bundle on an elliptic curve is naturally realized as a hyperplane in the dual of a Lie algebra, which comes with Kirillov-Kostant bracket.

In section four we use the language of $q$-difference modules to link the cohomology of $G$-bundles on $E_q$ with the corresponding multipliers in $LG_q$. Finally in section five we find an explicit formula for our bracket in the case of $\SL_2$.  We also describe how the symplectic leaves sit inside the moduli space.

Section six looks at loop groups. It was observed by Looijenga that twisted conjugacy classes in a loop group correspond to principal bundles on an elliptic curve, and we use this observation to provide another construction of the moduli space of parabolic bundles using a non-abelian Hamiltonian reduction. Using the theory of Poisson-Lie groups we construct a non-canonical Poisson bracket on the loop group compatible with the twisted conjugation action. For $\SL_2$ we are able to find an $r$-matrix so that the reduced Poisson structure is the same as that constructed in section three. We hope this can be generalized to other groups.

\section{Bundles on elliptic curves}
	\label{s:ellbun}

Vector bundles on an elliptic curve were classified by Atiyah in \cite{atiyah:57a}. The moduli space of principal bundles on an elliptic curve is well understood, see for example \cite{fm:98}. We are interested in the moduli space of bundles on an elliptic curve with additional parabolic structure. By this we mean a reduction of the structure group from $G$ to a parabolic subgroup $P\subset G$. In the case of $\GL_n$ a parabolic bundle is equivalent to a vector bundle along with a flag of sub-bundles. Thus the spaces we work with can be thought of as a kind of global Grassmannian defined over an elliptic curve.

We set the stage by recalling some of results from the appendix of Friedman-Morgan \cite{fm:02} where they study the moduli space of parabolic bundles in detail. Many of their results hold more generally than stated below, but for convenience we state all the assumptions we will be using up front.

Fix a complex reductive group $G$ with parabolic subgroup $P$. Let $U$ be the associated maximal unipotent subgroup of $P$ and $\pi:P\rightarrow P/U=L$ the Levi quotient. Note that $P$ is isomorphic to a semi-direct product of $L$ and $U$. 

Fix an $L$-bundle $\xi_0$ and consider the set of isomorphism classes of pairs $(\xi,\phi)$, where $\xi$ is a principal $P$-bundle on the elliptic curve $E$ and $\phi$ is an isomorphism $\phi:\xi/U\rightarrow \xi_0$. The isomorphism classes are classified by the non-abelian cohomology set $H^1(E,U(\xi_0))$, where $U(\xi_0)$ is the sheaf of unipotent groups $\xi_0\times_L U$.

Since $U$ is not abelian, {\it a priori} $H^1(E,U(\xi_0))$ is only a set, but in the specific situation we are considering, it can be shown to be an affine space. We are interested in the moduli space of parabolic bundles. Define moduli functor ${\bf F}$ as follows. Given a commutative $\bC$-algebra $S$, we have product space $X\times\Spec S$ with projections $\pi_1$ and $\pi_2$ to $X$ and $\Spec S$ respectively. Let ${\bf F}(S)$ be the set of isomorphism classes of pairs $(\Xi,\Phi)$ where $\Xi$ is a principal $P$-bundle over $X\times\Spec S$ and $\Phi:\Xi/U\rightarrow\pi^*_1\xi_0$ is an isomorphism. Thus
\begin{equation*}
	{\bf F}(S)=H^1(X\times\Spec S,U(\pi_1^*\xi_0)).
\end{equation*}
In general this functor is not representable, but in the case we are interested in the following theorem applies.

\begin{thm}
	\label{t:affmod}
	\cite{fm:02}
	Let $X$ be a projective scheme with $L$-bundle $\xi_0$, and let $U(\xi_0)$ be the associated sheaf of unipotent groups. Let $\{U_i\}_{i=1}^N$ be a decreasing filtration of $U$ by normal $L$-invariant subgroups such that for all $i$ the subquotient $U_i/U_{i+1}$ is contained in the center of $U/U_{i+1}$. Suppose that for all $i$
	\begin{equation*}
		H^0(X,(U_i/U_{i+1})(\xi_0))=0=H^2(X,(U_i/U_{i+1})(\xi_0)).
	\end{equation*}
Then
\renewcommand{\labelenumi}{\alph{enumi})}
\begin{enumerate}
	\item
	the cohomology set $H^1(X,U(\xi_0))$ has the structure of affine $n$-space $\bA^n$. More precisely there is a $P$-bundle $\Xi_0$ over $X\times\bA^n$ and an isomorphism $\Phi_0:\Xi_0/U\rightarrow\pi^*_1\xi_0$ such that the pair $(\Xi_0,\Phi_0)$ represents the functor ${\bf F}$ defined above.
	\item
	there is a natural action of the algebraic group $\Aut_L\xi_0$ on $H^1(X,U(\xi_0))$. This action lifts to an action on $\Xi_0$.
\end{enumerate}
\end{thm}

Since $P$ is a parabolic the lower central series of $U$ provides such a decreasing filtration. Higher cohomology vanishes on a curve so the only condition we need to impose on $\xi_0$ is that $H^0(E,(U_i/U_{i+1})(\xi_0))$ vanishes for all $i$.

\begin{rem}
	There is a marked point $0\in H^1(E,U(\xi_0))$ given by $(\xi_0\times_L P,I)$ where $I$ is the canonical identification of $(\xi_0\times_L P)/U$ with $\xi_0$. 
\end{rem}

Given an algebraic group $G$, let $\underline{G}$ denote the sheaf of morphisms from $E$ to $G$. Also, supposing that $\xi$ is a $P$-bundle lifting $\xi_0$, let $U(\xi)$ be the sheaf of sections of $\xi\times_P U\rightarrow E$. Then

\begin{lem}
	\label{l:moduni}
	\cite{fm:02}
	The cohomology $H^1(E,U(\xi))$ gives the set of all isomorphism classes of pairs $(\xi,\phi)$ where $\xi$ is a principal $P$-bundle and $\phi$ is an isomorphism $\phi:\xi/U\rightarrow \xi_0$.
\end{lem}

\begin{lem}
	\label{l:modpar}
	\cite{fm:02}
	There is a natural map $H^1(E, {\underline P})\rightarrow H^1(E,{\underline L})$ induced by the projection $\pi:P\rightarrow L$. The fiber over $\xi_0\in H^1(E,{\underline L})$ is the set of $P$-bundles on $E$ lifting $\xi_0$; it is given by $H^1(E,U(\xi))/H^0(E,L(\xi_0))$.
\end{lem}

Thus the set of parabolic bundles on $E$ fibers over the moduli space of principal $L$-bundles, where the fibers are the quotient of an affine space by a group action. The action is not in general well-behaved, but in the sequel we will focus on the smooth locus of the moduli space of parabolic bundles, where interesting additional structures arise.

\section{Poisson structure on parabolic bundles}
	\label{s:poisson}

In this section we construct the Poisson bivector on the moduli space of parabolic bundles. Casimirs are investigated, and finally an alternative construction using Atiyah-Bott reduction is presented. The Atiyah-Bott construction is helpful in understanding why such a Poisson bracket should exist. The additional flag structure of the parabolic bundles is preserved under reduction by quotienting out by a smaller-than-usual group, so that the symplectic leaves are larger than the orbits being reduced, and the natural Kirillov-Kostant bracket is inherited by the quotient.

\subsection{The Poisson bivector}

Let $\cM(L)$ denote the smooth locus of the coarse moduli space $\cM(L)$ of principal $L$-bundles on $E$. This is well understood, for example when $L$ is simply connected $\cM(L)\approx (E\otimes_\bZ \Lambda_L)/W$ where $\Lambda_L$ is the weight lattice and $W$ the Weyl group. This is a weighted projective space, not necessarily smooth. Another extreme case is when $L$ is a torus, corresponding to $P$ a Borel subgroup of $G$. In this case $\cM(L)$ is a product of copies of $E$. 

By Theorem \ref{t:affmod}, for each principal $L$-bundle $\xi_0$ in $\cM(L)$ we have associated affine space $H^1(E,U(\xi_0))$. Let
\begin{equation*}
	\cM(U):=\big\{(\xi_0,\xi,\phi)\,\vline\, \xi_0\in \cM(L) \mbox{ and }
	(\xi,\phi)\in H^1(E,U(\xi_0))\big\}
\end{equation*}
so that we have affine fibration $\cM(U)\rightarrow \cM(L)$. Lemma \ref{l:modpar} states that the desired moduli space of parabolic bundles is then given fibrewise by the quotient of $H^1(E,U(\xi))$ by $H^0(E,L(\xi_0))$. The action is not free and so the quotient is not well behaved.

\begin{eg}
	\label{eg:sl2flag}
	We consider the case of $\SL_2$	with parabolic subgroup $P$ the lower triangular matrices and $U$ the strictly lower triangular matrices. The Levi $L$ is isomorphic to $\bC^*$, so start by picking a line bundle $\xi_0$. A parabolic bundle corresponds to a flag
\begin{equation*}
	0\rightarrow \xi_0^*\rightarrow V\rightarrow \xi_0\rightarrow 0
\end{equation*}
where $V$ is a rank two degree zero vector bundle on $E$. The sheaf of unipotent groups $U(\xi)$ is $(\xi^*_0)^{\otimes 2}$. Theorem \ref{t:affmod} requires $H^0(E,(\xi_0^*)^{\otimes 2})=0$, meaning that $\xi_0$ must have positive degree, which we denote by $k$. The moduli space $H^1(E,U(\xi))$ is the $2k$-dimensional vector space $\Ext^1(\xi_0,\xi_0^*)$ classifying extensions $(V,\phi)$ where $\phi:V/\xi_0^*\xrightarrow{\approx}\xi_0$.

There is a natural $\bC^*$-action on $\Ext^1(\xi_0,\xi_0^*)$ given by rescaling $\phi$. Quotienting this out gives the moduli space of parabolic bundles lifting $\xi_0$. This is not a variety, but if we delete the origin, which is badly behaved under the $\bC^*$-action, the quotient is the projective space $\bP^{2k-1}$. Note the origin is the trivial extension.

Thus forgetting the trivial extensions, the moduli space of parabolics breaks into connected components labeled by the degree of the corresponding Levi-bundle $\xi_0$. Each component is a fibration $\cM(P)_k\rightarrow \Jac^k(E)$ with fibers isomorphic to $\bP^{2k-1}$.
\end{eg}

In general to avoid singular points as in Example \ref{eg:sl2flag}, we let $\cM(P)$ denote the smooth locus of the moduli space of parabolic bundles. Thus $\cM(P)$ is a fibration over $\cM(L)$, where the fiber over $\xi_0\in\cM(L)$ is a subspace of the quotient $H^1(E,U(\xi_0))/H^0(E,L(\xi_0))$.

Associated to the inclusion $P\hookrightarrow G$ is exact sequence of Lie algebras
\begin{equation*}
	0\rightarrow \fp\rightarrow \fg\rightarrow\fu_-\rightarrow 0.
\end{equation*}
Here we are using the Killing form to identify the quotient $\fg/\fp$ with the unipotent algebra $\fu_-$ opposite to $\fu$. Given a $P$-bundle $\xi$ denote by $\fp(\xi)$ the sheaf of infinitesimal automorphisms of $\xi$ respecting the parabolic structure, and let $\fg(\xi)$ be the sheaf of infinitesimal automorphisms of the $G$-bundle induced by $\xi$ under the inclusion $P\hookrightarrow G$. Slightly abusing notation we let $\fu_-(\xi)$ be the quotient $\fg(\xi)/\fp(\xi)$.

This gives rise to the following diagram of sheaves on $E$
\begin{equation}
	\label{e:diamond}
\xymatrix{
&\fu(\xi)\ar[d]\ar[dr]^0\ar[dl]_\rho&\\
\fp(\xi)\ar[r]&\fg(\xi)\ar[r]\ar[d]&\fu_-(\xi)\\
 & \fp_-(\xi)\ar[ur]_\tau\\
}\end{equation}
where $\fu(\xi)$ is the sheaf of Lie algebras associated to $U(\xi)$ and again $\fp_-(\xi)$ is the quotient by abuse of notation. The diagonal maps are induced by the corresponding maps on Lie algebras.

Recall that $\cM(P)$ refers to the smooth locus of the moduli space of parabolic bundles on $E$. By standard deformation theoretic arguments the tangent space at $\xi$ is $H^1(E,\fp(\xi))$, and by Serre duality combined with the Killing form the cotangent space is $H^0(E,\fp_-(\xi))$. Define the Poisson map
$B^\#_\xi:\Omega_\xi\cM(P)\rightarrow T_\xi\cM(P)$ as given by either side of the commutative diagram
\begin{equation}
	\label{e:bivec}
\xymatrix{
&  H^0(E,\fp_-(\xi))\ar[dr]^\delta\ar[dl]_\tau\ar[dd]^{B_\xi^\#} \\
H^0(E,\fu_-(\xi))\ar[dr]_\delta & & H^1(E,\fu(\xi))\ar[dl]^\rho \\
& H^1(E,\fp(\xi)).
}\end{equation}
where $\rho$ and $\tau$ are as in (\ref{e:diamond}) and the $\delta$'s are coboundary maps. Alternatively we can write the bivector in terms of the complex
\begin{equation*}
	C(\xi)=\left[\fg(\xi)\xrightarrow{d}(\fg/\fp)(\xi)\right]
\end{equation*}
with $d$ the obvious map. The Poisson map $B_\xi^\#$ can alternatively be realized as the cohomology of the map between complexes:
\begin{equation*}
\xymatrix{
 [\fu(\xi)\ar[d]_{B^\#_0}\ar[r]^{-d^*} & \fg(\xi)]\ar[d]_{B^\#_1} & H^1(E,C^*[-1])\ar[d]_{B_\xi^\#} \\
[\fg(\xi)\ar[r]^d & (\fg/\fp)(\xi)] & H^1(E,C) 
}
\end{equation*}
with $B^\#_0$ negative inclusion so the diagram commutes. The map of complexes $(B^\#_0,B^\#_1)$ is homotopic to $(B^\#_0,B^\#_1)[-1]$ with homotopy given by $h:\fg\xrightarrow{-2\cdot}\fg$. It follows by skew-symmetry of the pairing $H^1(C)\otimes H^1(C^*[-1])\rightarrow H^2(C\otimes C^*[-1])$ that $B_\xi^\#$ is skew-symmetric.

\begin{thm}
	The skew-symmetric tensor $B^\#:\Omega\cM(P)\rightarrow T\cM(P)$ defines a Poisson structure.
\end{thm}

\begin{proof}
	The map $B^\#$ is defined fibrewise for each $\xi\in \cM(P)$ above. It can be shown to define a global holomorphic map by adapting the reasoning in Mukai \cite{muk:84}. Alternatively, it follows by Proposition \ref{p:ab_red} below.
	
	All that remains is to prove the Jacobi identity. Bottacin \cite{bottacin:95} has shown the Jacobi identity for a skew-symmetric tensor $B^\#$ is equivalent to
	\begin{equation}
		\label{e:jacobi}
		B^\#(\omega^1)\cdot\bra B^\#(\omega^2),\omega^3\ket-
		\bra\left[B^\#(\omega^1),B^\#(\omega^2)\right],\omega^3\ket
		+ \mbox{ cyclic permuations }=0
	\end{equation}
	for any three 1-forms on $\cM(P)$, where $\bra\cdot,\cdot\ket$ is the pairing of vector fields with 1-forms. For us it is convenient to use the definition:
	\begin{equation*}
		B_\xi^\#:H^0(E,\fp_-(\xi))\xrightarrow{\delta}H^1(E,\fu(\xi)) \xrightarrow{\rho}H^1(E,\fp(\xi)).
	\end{equation*}
	Let $\pi_1$ and $\pi_2$ denote the projections of $E \times \cM(P)$ to $E$ and $\cM(P)$ respectively. Consider the extension on $E\times\cM(P)$
	\begin{equation*}
		0\rightarrow\pi_1^*\fp(\xi)\rightarrow\sD^1_{E\times\cM(P)/E} (\xi) \rightarrow \pi^*_2T\cM(P)\rightarrow 0
	\end{equation*}
	and sections of $\sD^1(\xi)$ are first-order differential operators with scalar symbol, $\pi_1^*\sO_E$-linear, and preserving the parabolic structure. This sequence can be pushed forward to a long exact sequence on $\cM(P)$
	\begin{equation*}
		\cdots\rightarrow T\cM(P)\xrightarrow{\approx}R^1\pi_{2*}\pi_1^*\fp(\xi) \xrightarrow{0} R^1\pi_{2*}\sD^1(\xi)\rightarrow\cdots
	\end{equation*}
	where the isomorphism is given pointwise by the identification of the tangent space $T_\xi\cM(P)$ with $H^1(E,\fp(\xi))$. From this it follows that the image in $R^1\pi_{2*}\sD^1(\xi)$ of an element in $R^1\pi_{2*}\fp(\xi)$ is a coboundary. Choose a sufficiently fine cover $\cU=\{U_i\}$ of $\cM(P)$ near $\xi$ so that we can write the 1-cocycle $B^\#_\xi(\omega^k)$ as coboundary $D^k_j-D^k_i=:D^k_{ij}$ on the intersection $U_{ij}=U_i\cap U_j$. Considered locally on subsets $U_{ij}$ skew symmetry is then
	\begin{equation*}
		\bra B^\#_\xi(\omega^h),\omega^l_j\ket=\bra D^h_{ij},\omega^l_i \ket = -\bra D^l_{ij},\omega^h_i\ket=-\bra B_\xi^\#(\omega^l), \omega^h_i\ket.
	\end{equation*}
	Finally, putting this together (\ref{e:jacobi}) is locally given by
	\begin{align*}
	& D^1_{ij}\cdot<{B^\#}(\omega^2),\omega^3_i> - <[{B^\#}(\omega^1),{B^\#}(\omega^2)],\omega^3_i> + \mbox{ cyc perm}\\
	 = & <D^1_{ij}\cdot{B^\#}(\omega^2),\omega^3_i> + <{B^\#}(\omega^2),D^1_{ij}\cdot\omega^3_i> - \\
	& - <D^1_{ij}\cdot{B^\#}(\omega^2) - D^2_{ij}\cdot{B^\#}(\omega^1),\omega^3_i>+ \mbox{ c.p.}\\
	= & <D^2_{ij}\cdot{B^\#}(\omega^1),\omega^3_i> + <{B^\#}(\omega^2),D^1_{ij}\cdot\omega^3_i> + \mbox{ c.p.}\\
	= &\,\, 0 \,\,\mbox{ by skew symmetry.}
	\end{align*}
The Leibniz rule is used in the first equality and combined with skew-symmetry gives the Jacobi.
\end{proof}

\begin{rem}
	\label{r:deg}
	In the degenerate case where $P=G$ the Poisson structure is trivial.
\end{rem}

\subsection{Casimir maps}
\label{ss:casimir}
The rich group-theoretic structure of the moduli space of parabolics can be used to construct Casmir maps for the Poisson structure.

The projection $P\rightarrow P/U=L$ induces a map
\begin{align*}
	\det:\cM(P) & \rightarrow\cM(L) \\
	\xi & \mapsto \xi/U=:\xi_0
\end{align*}
to the moduli space of $L$-bundles on $E$. The differential $\det_*:T\cM(P)\rightarrow T\cM(L)$ and its dual are induced by maps of sheaves on the ``completion'' of diagram (\ref{e:diamond}):

\begin{equation}
\label{e:det}
\xymatrix{
& & \fu(\xi)\ar[d]\ar[dr]^0\ar[dl]_\rho & \\
& \fp(\xi)\ar[r]\ar[dl]_{\det_*} & \fg(\xi)\ar[r]\ar[d]&\fu_-(\xi) \\
\fl(\xi/U) & & \fp_-(\xi)\ar[ur]_\tau \\
& \fl(\xi/U)\ar[ur]_{\det^*}
}
\end{equation}

\begin{lem}
	\label{l:det}
	The determinant map $\det:\cM(P)\rightarrow \cM(L):\xi\mapsto\xi/U$ is a Casimir map, so the subvarieties $\cM(P)_{\xi_0}:=\det^{-1}(\xi_0)$ are Poisson subvarieties.
\end{lem}

\begin{proof}
	The following simple proof was suggested by the referee. Recall, diagram (\ref{e:bivec}), that $B^\#_\xi$ is the composition $\delta\circ\tau$. Composing $\fl(\xi_0)\rightarrow \fp_-(\xi)$ with $\tau:\fp_-(\xi)\rightarrow \fu_-(\xi)$ gives the zero map, implying $\det$ is a Casmir.
\end{proof}

We now investigate inclusions of parabolics. Suppose $P_1\subset P_2\subset G$. Then there is the forgetful map $f:\cM(P_1)\rightarrow \cM(P_2)$.

\begin{lem}
	\label{l:incl}
	The forgetful map $f:\cM(P_1)\rightarrow \cM(P_2)$ induced by inclusion is a Poisson map.
\end{lem}

\begin{proof}
	We have to show
	\begin{equation*}
		\xymatrix{
		H^0(\fg/\fu_1(\xi))\ar[r]^{B^\#_1} & H^1(\fp_1(\xi))\ar[d] \\
		H^0(\fg/\fu_2(\xi))\ar[u]\ar[r]^{B^\#_2} & H^1(\fp_2(\xi))
		}
	\end{equation*}
	commutes. This follows from commutative diagram
	\begin{equation*}
		\xymatrix{
		\fp_1\ar[d] & \fu_1\ar[l]\ar[d] & \fu_2\ar[l] \\
		\fp_2\ar[r] & \fg\ar[r]\ar[d] & \fg/\fp_2 \\
		\fg/\fu_2\ar[r] & \fg/\fu_1\ar[r] & \fg/\fp_1\ar[u]
		}
	\end{equation*}
	which reduces the problem to showing
	\begin{equation*}
		\xymatrix{
		H^0(\fg/\fp_1)\ar[r]\ar[d] & H^1(\fp_1)\ar[d] \\
		H^0(\fg/\fp_2)\ar[r] & H^1(\fp_2)
		}
	\end{equation*}
	commutes, which is clear.
\end{proof}

Finally there is forgetful map $f:\cM(P)\rightarrow\cM(G)$ taking a $P$-bundle to a $G$-bundle by forgetting the flag structure. Lemma \ref{l:incl} states this is a Poisson map, and since by Remark \ref{r:deg}
$\cM(G)$ has trivial Poisson structure $f$ is a Casimir and cuts out Poisson submanifolds.

The Casimirs can be depicted as follows
\begin{equation*}
	\label{e:casimirs}
	\xymatrix{
	& \cM(P)\ar[dl]_f\ar[dr]^\det \\
	\cM(G) & & \cM(L)
	}
\end{equation*}

\subsection{Atiyah-Bott reduction}

In this subsection we provide an alternate construction of the Poisson bivector. This point of view also motivates why we would expect a Poisson structure in the first place.

Atiyah-Bott reduction \cite{ab:83} is a method for constructing the moduli space of holomorphic $G$-bundles on a Riemann surface $M$. The space $\sC_G$ of holomorphic structures on a fixed principal bundle is an affine space modelled on the vector space $\Omega^{0,1}(M,\fg)$. The moduli space of holomorphic bundles is then a quotient by the action of the gauge group $\sG=C^\infty(M,G)$.

Suppose $G$ is simple and simply connected and our Riemann surface is an elliptic curve $E$. Then topologically there is only the trivial $G$-bundle $G\times E$ on $E$. Let $\sE G=C^\infty(E,G)$ be the gauge group with Lie algebra $\sE\fg=\Omega^0(E,\fg)$. Fix a holomorphic differential $\eta\in H^0(E,\omega_E)$ and following \cite{ef:94} define a central extension $\widehat{\sE\fg}$ using cocycle $\Omega_\eta(X,Y)=\int_E \eta\wedge\bra X, \dbar Y\ket$. Let $(\widehat{\sE\fg})^* = \{\lambda \dbar + \xi\,\vline\,\xi\in\Omega^{0,1}_E\otimes\fg\}$ be the dual, where we identify $\fg\approx\fg^*$ using the Killing form. We have pairing
\begin{align*}
	\widehat{\sE\fg}\otimes(\widehat{\sE\fg})^* & \rightarrow \bC \\
	(\mu k + X, \lambda\dbar + \xi) & \mapsto \int_E\eta\wedge\bra \xi, X\ket - \lambda\mu.
\end{align*}
The gauge groups acts on $(\widehat{\sE\fg})^*$ via a twisted co-adjoint action
\begin{equation*}
	h\circ(\lambda\dbar+\xi):= \lambda\dbar+h^{-1}(\dbar h)+ \Ad (h\circ\xi).
\end{equation*}
This actions leaves the hyperplane fixing $\lambda$ invariant and in particular setting $\lambda=1$ we recover $\sC_G=\{\dbar+\xi\,\vline\,\xi\in\Omega^{0,1}_E\otimes\fg\}$, the infinite dimensional affine space of operators $D=\dbar+\xi$ on $E\times G$. 

It is a classical theorem of Atiyah and Bott that the orbits of $\sE G$ on $\sC_G$ are in bijective correspondence with holomorphic principal $G$-bundles on $E$. The case of an elliptic curve is special since  $\sC_G$ identifies naturally with an affine hyperplane in the dual of a Lie algebra. Therefore it comes with a natural Kirillov-Kostant Poisson bivector:

\begin{align}
	\label{e:dlpoisson}
	\widehat{\sE\fg} & \xrightarrow{B^\#_{\dbar+\xi}} (\widehat{\sE\fg})^* \\
	(\mu k + X) & \mapsto \left((\nu k + Y) \mapsto \int_E\eta\wedge\big\langle\dbar X+[\xi,X],Y\big\rangle\right).
\end{align}
This restricts to a Poisson bivector $\sE\fg\rightarrow \sC_G$. Symplectic leaves are cut out by co-adjoint orbits of the gauge group $\sE G$, and thus by the result of Atiyah and Bott are in bijection with $G$-bundles on $E$.

Now consider diagram
\begin{equation*}
	\xymatrix{
	\sC_L\ar[d]^{\mu_L} & \sC_G\ar[d]^{\mu_P} \\
	\Omega^{0,1}(E,\fl^*)\ar[r] & \Omega^{0,1}(E,\fp^*)
	}
\end{equation*}
where the vertical arrows are moment maps associated to gauge actions of $\sE L$ and $\sE P$, and the horizontal map is induced by the projection $\pi:P\rightarrow L$. Moment map $\mu_L$ goes $\dbar +\xi\mapsto \xi$, whereas $\mu_P$ is $\dbar +\xi\mapsto \xi$ composed with the projection induced by $\fg^*\rightarrow \fp^*$.

The moduli space of parabolic bundles is constructed as follows. Pick a $L$-bundle $\xi_0$. This corresponds to a $\sE L$-orbit $\cO_{\xi_0}$ in $\sC_L$ which projects into $\Omega^{0,1}(E,\fl)$ and then can be included into $\Omega^{0,1}(E,\fp^*)$. Since $P=L\cdot U$ is a semidirect product $\mu_L(\cO_{\xi_0})\subset \Omega^{0,1}(E,\fp^*)$ is invariant under the action of gauge group $\sE P$. It follows by the theory of Atiyah-Bott that the $\sE P$ orbits in $\mu_P^{-1}(\cO_{\xi_0})$ are in bijection with parabolic bundles. Again we have problems related to bad-quotients, but in specific examples we can restrict to orbits where the quotient is well-behaved, see Example \ref{eg:sl2flag}.

\begin{prop}
	\label{p:ab_red}
	Atiyah-Bott reduction takes the Kirillov-Kostant Poisson bivector (\ref{e:dlpoisson}) on $\sC_G$ to the parabolic Poisson bivector $B^\#$ (\ref{e:bivec}) on $\cM(P)_{\xi_0}$.
\end{prop}

\begin{proof}
	The bivector (\ref{e:bivec}) is constructed on cohomology groups using the coboundary operator $\delta$ and map $\tau$. The Kirillov-Kostant bivector is constructed using the Dolbeault operator $\dbar+\xi$ and the idea of the proof is to show that these coincide under the Atiyah-Bott reduction. The leaves of the Kirillov-Kostant bracket are $\sE G$-orbits, so reducing by the smaller group $\sE P$ preserves the Poisson structure.
	
	Recall we identified $\fp^*$ with $\fp_-$ and $\fu^*$ with $\fu_-$ using the Killing form. The parabolic bivector (\ref{e:bivec}) is composition
	\begin{equation*}
		H^0(E,\fp^*(\xi))\xrightarrow{\tau} H^0(E,\fu^*(\xi)) \xrightarrow{\delta} H^1(E,\fp(\xi)).
	\end{equation*}
	This can be expressed in terms of Dolbeault cohomology. Map $\tau$ is induced by inclusion $\fu\hookrightarrow \fp$, and the differential $\delta$ is given by $\dbar +\tilde{\xi}$ which defines the complex structure corresponding to $G$-bundle $\xi$, where here $\tilde{\xi}\in \Omega^{0,1}(E,\fp)$.
\end{proof}

\begin{cor}
	\label{c:leaves}
	Symplectic leaves in $\cM(P)_{\xi_0}$ are given by fibers of the forgetful map $f:\cM(P)_{\xi_0}\rightarrow \cM(G)$.
\end{cor}

\section{Difference modules}
	\label{s:diff}

We wish to have a more explicit description of the Poisson bracket constructed above. For this purpose we turn to the language of $q$-difference modules which provides a means for explicitly writing out sections of vector bundles on an elliptic curve $E_q$. Our basic reference is \cite{vdpr:qa05}.

The connection between $q$-difference modules and elliptic curves is that any elliptic curve $E$ can be written as a quotient $\bC^*/q^\bZ$ for some $q$ with norm less than one. Vector bundles on $E$ correspond naturally to $q$-difference modules and their cohomology can be naturally calculated. The multiplicative quotient $\bC^*/q^\bZ$ also directly leads from elliptic curves to conjugacy classes in loop groups.

Thus $q$-difference modules provide a bridge between the sheaf-theoretic construction of the Poisson bracket in the previous section, and the group-theoretic construction to be presented in \S\ref{s:loop}.

\subsection{Twisted conjugacy classes}
	\label{ss:qconj}

An elliptic curve is often written as a quotient of $\bC$ by a rank two lattice $\bZ\bra1,\tau\ket$ where $\tau$ lies in the upper half plane. For our purposes it is useful to replace this additive point of view with a multiplicative one: start with map $\bC\rightarrow \bC/\bZ=\bC^*$ given by $e^{2\pi i-}$. This sends $\tau$ to a non-zero $q$ with $|q|<1$. The elliptic curve is given by projection $\bC^*\rightarrow \bC^*/q^\bZ$ and is denoted by $E_q$. Given $\eta\in\bC$ write $\bar{\eta}\in E_q$ for the its image under projection.

Given complex Lie group $G$, write $LG$ for the loop group of holomorphic maps from $\bC^*$ to $G$. This group has natural $\bC^*$-action given by twisting or rotation: $q\in\bC^*$ acts by
\begin{equation}
	q:a(z)\mapsto a_q(z):=a(qz).
\end{equation}
Write $\widetilde{LG}:=\bC^*\ltimes LG$ for the semidirect product. The (right) conjugation action of the subgroup $LG$ on $\widetilde{LG}$ goes as follows
\begin{equation}
	\label{e:qconj}
	(1,g(z))^{-1}\cdot(q,a(z))\cdot(1,g(z)
	= \left(q,g^{-1}(qz)\cdot a(z)\cdot g(z)\right).
\end{equation}
We fix a choice of $q\in\bC^*$ with $0\neq|q|<1$ and denote $\{q\}\times LG$ by $LG_q$. This is isomorphic to $LG$ and carries the twisted conjugation action by the loop group: $g:a\mapsto {^g}\!a=g^{-1}(qz)\cdot a(z)\cdot g(z)$. We will refer to the orbits of this action as $q$-conjugacy classes in $LG$. These have a surprising geometric significance discovered by Looijenga

\begin{thm}
	\label{t:looijenga}
	Let $G$ be a connected complex Lie group. There is a natural bijection between the set of $q$-conjugacy classes in $LG$ and the set of isomorphism classes of holomorphic principal $G$-bundles on $E_q$ with left $G$-action.
\end{thm}

\begin{proof}
	\cite{bg:96}
	The $q$-conjugation action of $LG$ on itself is given by (\ref{e:qconj}). Given an element $a(z)\in LG$ construct a principal $G$-bundle $\eta(a)$ as follows. Define a $q^\bZ$ equivariant structure on the trivial bundle $G\times\bC^*\rightarrow\bC^*$ by setting 
	\begin{equation}
		\label{e:qeqv}
		q:(g(z),z)\mapsto(g_q(z) a(z),q z).
	\end{equation}
	 This descends to give a principal $G$-bundle on $E_q$. The element $a(z)$ is known as the {\em multiplier} associated to $\eta(a)$; $q$-conjugate multipliers give rise to isomorphic $G$-bundles.

	Conversely, any principal $G$-bundle on $\bC^*$ is holomorphically trivial so principal $G$-bundles on $E_q$ pull back to trivial principal bundles on $\bC^*$ with $q^\bZ$-action. Fixing a trivialization this is of the form above. Changing trivializations corresponds to $q$-conjugation.
	
	Finally there is a natural left $G$-action on $G\times\bC^*$ given by $h:(g(z),z)\mapsto (h\cdot g(z),z)$ which descends to a left $G$-action on $\eta(a)$.
\end{proof}

\begin{rem}
	We use a {\it right} action for $q$-conjugation $LG_q\times LG\rightarrow LG_q$ so the multiplier acts in a slightly different way than in \cite{bg:96}.
\end{rem}

As an immediate consequence we observe

\begin{cor}
Given $a(z)\in LG$ let $Z_q(a)=\{g\in LG\,|\, g^{-1}_q\cdot a\cdot g=a\}$ be its centralizer in $LG$ under $q$-conjugation. There is a natural isomorphism $\Aut_G\eta(a)\rightarrow Z_q(a)$.
\end{cor}

\begin{proof}
An automorphism of $\eta(a)$ is a $G$-equivariant isomorphism $g:\eta\rightarrow\eta$ over the base $E_q$. Lifted to $\bC^*$ this is a map $\bC^*\rightarrow G$ which $q$-commutes with $a$.
\end{proof}

\subsection{An equivalence of categories}

Given an element $a(z)$ of $LG_q$ we have seen that we can construct a left $G$-bundle on $E_q$ with $a(z)$ as multiplier. Given a $G$-bundle $\eta$ and a right representation $\rho:G\rightarrow \GL(V)$ we can construct vector bundle $V\times_\rho \eta$ on $E_q$. The multiplier can be used to define a $q$-difference module structure on the space of holomorphic loops in $V$, and we use this point of view to understand the cohomology of $V\times_\rho\eta$.

A difference ring is a commutative ring $R$ with given automorphism $\phi$. A difference module is then a (left) free $R[\phi,\phi^{-1}]$ module of finite rank. In our case let $\sR$ be holomorphic functions on $\bC^*$ and $\phi$ be the automorphism $f(z)\mapsto f(qz)=f_q(z)$. A $q$-difference module is by definition a difference module over $\sR[\phi,\phi^{-1}]$.

Given a right representation $\rho:G\rightarrow GL(V)$, we construct a $q$-difference module for each $a\in LG$ as follows. Let $\sV=V\otimes_\bC\sR$ be the space of holomorphic function on $\bC^*$ with values in $V$. Already this has a natural $q$-difference module structure given $\phi\cdot{\bf f}(z)={\bf f}(qz)$. This corresponds to $1\in LG$. Given multiplier $a\in LG$ we twist this
\begin{equation}
	\label{e:multtwist}
	\phi_a\cdot {\bf f}(z):=(\phi {\bf f})(z)\cdot a={\bf f}(qz)\cdot a.
\end{equation}
Denote the new $q$-difference module by $\sV(a)$, and refer to $a$ as the multiplier of $\sV$. In examples $G$ will be $\SL_n$ and the representation will be the fundamental representation, with vectors written as rows and the group acting by multiplication on the right.

\begin{defn}
	Given a $q$-difference module $\sV$ we can naturally associate a vector bundle $\tilde{V}$ on $E_q$. Its sheaf of sections is
	\begin{equation*}
		\tilde{V}(U) = \{f\in\sO(\pi^{-1}U) \otimes_{\sR}\sV\, |\,\phi_\sV(f)=f\},
	\end{equation*}
	where $\pi:\bC^*\rightarrow E_q$ is the projection and given $U'\subset \bC^*$, $\sO(U')$ is the algebra of holomorphic functions on $U'$. This extends naturally to a functor $\sF$ from the category of finitely generated  $q$-difference modules to the category of vector bundles on $E_q$.
\end{defn}

\begin{rem}
	\cite{vdpr:qa05}
The functor $\sF$ is an equivalence of categories and extends to an equivalence between finitely generated not necessarily free $q$-difference modules and coherent sheaves on $E_q$.
\end{rem}

\begin{lem}
	Given $V$ a representation of $G$ and multiplier $a\in LG$, the vector bundles $\sF(\sV(a))$ and $V\times_G\eta(a)$ are naturally isomorphic.
\end{lem}

\begin{proof}
	The $q$-difference module action is given by (\ref{e:multtwist}). For the principal $G$-bundle $\eta(a)$, recall the $q$-equivariant structure (\ref{e:qeqv})
	\begin{equation*}
		(g(z),z)\mapsto (g_q(z)a(z),qz).
	\end{equation*}
	where $g(z)$ is a trivialization of $\eta(a)$ lifted to $\bC^*\times G$. Pick (trivial) trivialization $g(z)=1$; other trivializations give the same bundle, but with different multiplier. Local sections of the bundle $V\times_G \eta(a)$ pulled back to $\bC^*\times V$ are elements of $\sO(\pi^{-1}U)\otimes_\sR \sV$. The $q$-twist is given by $f(z)\mapsto f(qz)\cdot 1\cdot a$.
\end{proof}

\begin{lem}
\cite{vdpr:qa05}
There are isomorphisms $\ker(\phi_\sV-1)\rightarrow H^0(E_q,\tilde{V})$ and $\cok(\phi_\sV-1)\rightarrow H^1(E_q,\tilde{V})$.
\end{lem}

\begin{proof}
The kernel and cokernel fit into exact sequence
\begin{equation}
\label{e:qcohom}
	0\rightarrow \ker\rightarrow \sV\xrightarrow{\phi-1}\sV\rightarrow\cok\rightarrow 0
\end{equation}
where $(\phi-1)$ is the coboundary operator, see \cite{vdpr:qa05}.
\end{proof}

\begin{lem}
	\label{l:dualmult}
Given $q$-difference module $\sV(a)$ and picking some basis, the dual $q$-difference module $\sV(a)^*$ is given by $\sV(a^*)$ where $a^*=\bar{a}^{-1}$ is the inverse transpose of $a$.
\end{lem}

\begin{rem}
\label{r:loopserre}
Serre duality can be understood from this point of view. Exact sequence
\begin{gather*}
	0\rightarrow\ker(\phi_a-1)\rightarrow \sV(a)\xrightarrow{\phi_a-1}\sV(a)\rightarrow \cok(\phi_a-1)\rightarrow 0\,\,\mbox{ has dual}\\
	0\rightarrow\cok(\phi_a-1)^*\rightarrow\sV(a)^*\xrightarrow{(\phi_a-1)^*}\sV(a)^*\rightarrow\ker(\phi_a-1)^*\rightarrow 0\,\,\mbox{ which identifies with}\\
	0\rightarrow \ker(\phi_{a^*}-1)\rightarrow\sV(a^*)\rightarrow \sV(a^*)\rightarrow \cok(\phi_{a^*}-1)\rightarrow0.
\end{gather*}
Thus $H^0(E_q,\eta(a))^*\approx H^1(E_q,\eta(a^*))$ and similarly for $H^1(E_q,\eta(a))$.
\end{rem}

\subsection{Sections of line bundles}

We use the above framework to investigate sections of line bundles on $E_q$. Under equivalence $\sF$ line bundles correspond to rank one $q$-difference modules. These are of the form $L(\eta z)^k=\sR\cdot e$ with $\phi\cdot e = \eta^k z^k$ for some $\eta\in\bC^*$ and $k\in\bZ$. Sections of $\tilde{L}(\eta z)^k$ correspond to solutions of the functional equation in $\sR$
\begin{equation*}
	\eta^kz^kf_q-f=0
\end{equation*}
The trivial line bundle is given by $\tilde{L}(1)$ and sections are solutions of $f_q=f$, the constant functions.

The next simplest example is the line bundle $\tilde{L}(z)$. The equation $zf_q-f=0$ has unique solution up to scalar, the Jacobi theta function $\vartheta(z)=\sum_{l\in\bZ}q^{l(l-1)/2}z^l$. Sections of other line bundles of degree one $\tilde{L}(\eta z)$ are given by solutions of $\eta z f_q-f=0$. These are scalar multiples of $\vartheta(z,\eta)=\vartheta(\eta z)=\sum_{l\in\bZ}q^{l(l-1)/2}\eta^l z^l$.

Sections of line bundles of higher degree can also be described in terms of power series. Sections of $\tilde{L}(\eta z)^k$ are solutions of $\eta^k z^k f_q-f=0$. The space of solutions is a $k$-dimensional vector space. One choice of basis, which we will use in the sequel, is the following
\begin{equation}
	\label{e:theta}
	\vartheta_n^k(z,\eta)=\sum_{l\in\bZ}q^{ln}q^{kl(l-1)/2}\eta^{kl}z^{kl+n}
	\,\,\,\mbox{ for } n = 0,\ldots,k-1.
\end{equation}

\subsection{Degree one cohomology and Serre duality}

Line bundles with negative degree have non-zero first cohomology. The image of the operator $\phi-1$ on $L(\eta z)^{-k}$ contains functions of the form $\frac{f_q}{\eta^k z^k}-f$. Serre duality implies the cokernel must be a $k$-dimensional vector space, and it is easy to check that a basis of representatives is given by $\{[1],\ldots,[z^{k-1}]\}$, with relations
\begin{equation*}
	\eta^{kl}[z^{kl+n}]=q^{ln}q^{kl(l+1)/2}[z^n].
\end{equation*}

Given a line bundle $L$ on $E_q$ by Serre duality there is a pairing $H^0(E_q,L)\otimes H^1(E_q,L^*)\rightarrow H^1(E_q, \sO)=\bC$.

\begin{lem}
In terms of rank one $q$-difference modules the Serre pairing is given by
\begin{align*}
	H^0(E_q,\tilde{L}(\eta z)^k)\otimes H^1(E_q,\tilde{L}(\eta z)^{-k}) & \rightarrow H^1(E_q,\tilde{L}(1))\rightarrow \bC \\
	\vartheta^k(z,\eta)\otimes [f(z)] &\mapsto [\vartheta\cdot f] \mapsto (\vartheta\cdot f)(0).
\end{align*}
\end{lem}

\begin{proof}
It is clear that $\tilde{L}(\eta z)^k$ and $\tilde{L}(\eta z)^{-k}$ are dual, see also Remark \ref{r:loopserre}. To check the pairing is well-defined, match element
\begin{equation*}
	\vartheta^k_j(z,\eta)=\sum_l q^{jl}q^{kl(l-1)/2}\eta^{kl}z^{kl+j}\,\,\mbox{ with }\,[q^{k\alpha+k-j}z^{k\alpha-j}-\eta^k z^{k(\alpha+1)-j}]=0.
\end{equation*}
We only consider the pairing with terms of this form since replacing the $-j$ with something else modulo $k$ will give zero when we consider the constant term. Multiplying out gives
\begin{equation*}
	\sum_l\left(q^{k\alpha+k+j(l-1)}q^{kl(l-1)/2}\eta^{kl}z^{k(l+\alpha)}-q^{jl}q^{kl(l-1)/2}\eta^{k(l+1)}z^{k(l+\alpha+1)}\right)
\end{equation*}
so substituting $l=-\alpha$ and $l=-\alpha-1$ respectively restricts attention to the constant term and we obtain
\begin{equation*}
	q^{k\alpha+k-j(\alpha+1)}q^{k\alpha(\alpha+1)/2}\eta^{-k\alpha}-q^{j(-\alpha-1)}q^{k(\alpha+1)(\alpha+2)/2}\eta^{-k\alpha}=0.
\end{equation*}
\end{proof}

\begin{rem}
	\label{r:functionals}
Serre duality gives us another way to think about $H^0(E_q,\tilde{L}(\eta z)^k)$. The group $H^1(E_q,\tilde{L}(\eta z)^{-k})$ is given by functions in $\sR$ modulo the relation $[q^lz^{l-k}]=[\eta^k z^l]$. Thus elements of $H^0(E_q,\tilde{L}(\eta z)^k) = H^1(E_q,\tilde{L}(\eta z)^{-k})^*$ should be linear functionals on $\sR$ invariant under this relation.

Direct computation finds a basis for these of the form
\begin{equation}
\label{e:thetafunctional}
\theta^k_n(\eta^{-1})=\sum_{l\in\bZ}q^{nl}q^{kl(l+1)/2}\eta^{-kl}\varphi_{kl+n}
\end{equation}
where $\varphi_n(f)$ is the $n^{th}$ coefficient of the power series expansion of $f$. Note these functionals are closely related to the functionals arising when working with functions on the loop group invariant under $q$-conjugation, Example \ref{e:loopfunctional}. Finally, evaluating against the basis $\{[1],\ldots,[z^{k-1}]\}$ of $H^1(E_q,\tilde{L}(\eta z)^{-k})$ we observe
\begin{equation}
\label{e:evaltheta}
\bra \vartheta^k_{-n},[z^m]\ket=\delta_{mn}=\bra\theta^k_n,[z^m]\ket.
\end{equation}
\end{rem}

\section{Rank two vector bundles}
	\label{s:ranktwo}

The moduli space of rank two vector bundles with a flag structure are the simplest examples where the Poisson bracket can be explicitly calculated and is non-trivial. For this reason we investigate these bundles in detail. Using the language of $q$-difference modules we analyze the structure of extensions of line bundles on an elliptic curve, and apply our results to write down the Poisson bracket. Finally we describe the symplectic leaves of the bracket in this simplest case.

\subsection{Extensions and difference modules}

For $\SL_2$ there is up to conjugation only one choice of parabolic, so we fix $P$ as the lower triangular matrices. Then $U$ is the abelian group of strictly lower triangular matrices and the Levi is $\bC^*$. Recall from Example \ref{eg:sl2flag} that the moduli space $H^1(E_q,U(\xi))$ classifies extensions
\begin{equation}
	\label{e:rk2flag}
	0\rightarrow \xi^*_0\rightarrow V\rightarrow\xi_0\rightarrow 0
\end{equation}
and so is isomorphic to $\Ext^1_{E_q}(\xi_0,\xi_0^*)$. Forgetting the isomorphism $V/\xi_0^*\rightarrow\xi_0$ corresponds to quotienting out the obvious $\bC^*$ action. Under the functorial equivalence between vector bundles $V$ and $q$-difference modules $\sV$ this extension corresponds to multiplier
\begin{equation*}
	a=\left(\begin{matrix}
		(\eta z)^{-k} & 0 \\
		(\eta z)^k x & (\eta z)^k
	\end{matrix}
	\right)
	\,\mbox{ acting on}\,
	\left(\begin{matrix}
		e_1 & e_2
	\end{matrix}
	\right),
\end{equation*}
where recall (\ref{e:multtwist}) the associated $q$-difference module has difference action $\phi_a\cdot {\bf f}(z)={\bf f}(qz)\cdot a$. 
Explicitly, the action of $\phi$ on $\sV$ is
\begin{equation*}
	\phi:f\cdot e_1+g\cdot e_2\mapsto
	\left(\frac{f_q}{\eta^k z^k}+ \eta^k z^k x g_q
	\right)\cdot e_1
	+ (\eta z)^kg_q\cdot e_2.
\end{equation*}
The off-diagonal entry represents the extension class, see below, normalized to produce nicer formulae. We prove two lemmas calculating coboundary maps and extensions classes coming from (\ref{e:rk2flag}).

\begin{lem}
	\label{l:cobdry}
	The coboundary map in the long exact sequence coming from (\ref{e:rk2flag}) is
	\begin{align*}
		\delta:H^0(E_q,\xi_0) & \rightarrow H^1(E_q,\xi_0^*) \\
		b & \mapsto [b x]
	\end{align*}
\end{lem}

\begin{proof}
	The recipe is as follows: a section $b$ of $\xi_0=\tilde{L}(\eta z)^k$ satisfies $\eta^k z^k b_q=b$. Lift it to an element of $\sV$, $\tilde{b}=\tilde{a}\cdot e_1+b\cdot e_2$ for arbitrary choice of $\tilde{a}\in\sR$. Now consider
	\begin{equation*}
		(\phi-1)(\tilde{b})=
		\left(\frac{\tilde{a}_q}{\eta^k z^k} -\tilde{a} + \eta^k z^k xb_q
		\right)\cdot e_1+
		\left(\eta^k z^k b_q-b
		\right)\cdot e_2=\left(
		\frac{\tilde{a}_q}{\eta^k z^k}-\tilde{a}+\eta^k z^k xb_q
		\right)\cdot e_1.
	\end{equation*}
	This is a representative of the cohomology class $[bx]$.
\end{proof}

\begin{lem}
	\label{l:extclass}
	The extension (\ref{e:rk2flag}) is represented by cohomology class
	$[x]\in \Ext^1_{E_q}(\xi_0,\xi_0^*)$.
\end{lem}

\begin{proof}
	Consider the long exact sequence associated to $\Hom_{E_q}(\xi_0,-)$:
	\begin{equation*}
		\cdots\rightarrow\Hom(\xi_0,V)\rightarrow\Hom(\xi_0,\xi_0)
		\xrightarrow{\delta}\Ext^1(\xi_0,\xi_0^*)\rightarrow\cdots
	\end{equation*}
	The extension class is given by $\delta(Id)$. Tensoring (\ref{e:rk2flag}) with $\xi_0^*$ gives
	\begin{equation}
		\label{e:sesext}
		0\rightarrow (\xi_0^*)^{\otimes 2}\rightarrow \xi_0^*\otimes V
		\rightarrow \sO\rightarrow0.
	\end{equation}
	The multiplier for the extension is
	\begin{equation*}
		\left(\begin{matrix}
			(\eta z)^{-2k} & 0 \\
			x & 1
		\end{matrix}
		\right)
	\end{equation*}
	so applying the recipe from Lemma \ref{l:cobdry} gives $\delta(Id)=[x]$.
\end{proof}

To describe the Poisson bivector we need to investigate infinitesimal automorphisms preserving the flag (\ref{e:rk2flag}). First by Lemma \ref{l:dualmult} the dual bundle $V^*$ has multiplier
\begin{equation*}
	\left(\begin{matrix}
		(\eta z)^k & -(\eta z)^k x \\
		0 & (\eta z)^{-k}
	\end{matrix}
	\right),
\end{equation*}
so the rank four vector bundle $\underline{\End}(V)$ has multiplier
\begin{equation*}
	\left(\begin{matrix}
		1 & 0 & -x & 0 \\
		(\eta z)^{2k}x & (\eta z)^{2k} & -(\eta z)^{2k}x^2 &-(\eta z)^{2k} x \\
		0 & 0 & (\eta z)^{-2k} & 0 \\
		0 & 0 & x & 1
	\end{matrix}
	\right)
	\,\,\mbox{ acting on }\,\,
	\left(\begin{matrix}
		e_1^*\otimes e_1 \\
		e_1^*\otimes e_2 \\
		e_2^*\otimes e_1 \\
		e_2^*\otimes e_2
	\end{matrix}
	\right)^T.
\end{equation*}
From this we read off the sections of $\underline{\End}(V)$ satisfy
\begin{align*}
	a_q-a & = -bx &  (\eta z)^{2k}b_q & = b \\
	\frac{c_q}{(\eta z)^{2k}}-c & = ax-dx-bx^2 & d_q-d & = bx
\end{align*}

Recall exact sequence of sheaves from diagram (\ref{e:det})
\begin{equation}
	\label{e:ses}
	0\rightarrow\fu(\xi)\rightarrow\fp(\xi) \rightarrow \fl(\xi_0) \rightarrow 0.
\end{equation}
In this situation $\fu(\xi)=\underline{\Hom}_{E_q}(\xi_0,\xi^*_0)$ and $\fl(\xi_0)=\underline{\End}_{E_q}(\xi_0)$. Let $\underline{\End}_0(V)$ be the sheaf of sections of $\underline{\End}(V)$ with trace zero. Sheaf $\fp(\xi)\subset \underline{\End}_0(V)$ consists of flag-preserving infinitesimal automorphisms, and is obtained by forcing $b=0$. Thus extension (\ref{e:sesext}) is the instantiation of (\ref{e:ses}) for $\fg=\fsl_2$.

\begin{prop}
	\label{p:qbivec}
	The Poisson bivector (\ref{e:bivec}) is the composition
	\begin{equation*}
	\begin{matrix}
	B^\#:H^0(E,\fp_-) & \xrightarrow{\tau} & H^0(E,\fu_-) & \xrightarrow{\delta} & H^1(E,\fp) \\
	\\
	\left(\begin{matrix}
	a & b' \\
	0 & -a\end{matrix}\right) & \mapsto  & \left(\begin{matrix}
	0 & b'\\
	0 & 0\end{matrix}\right) & \mapsto & \left[\begin{matrix}
	0 & 0 \\
	(b'x-2a)x & 0
	\end{matrix}\right]
	\end{matrix}
	\end{equation*}
\end{prop}

\begin{proof}
	Note the bivector is zero for the trivial extension $x=0$. The sheaf $\fp_-(\xi)$ is given by forcing $c=0$. Then the equations above simplify to
	\begin{equation*}
		\begin{matrix}
			a_q-a=-bx, & &
			d_q-d=bx, & & 
			(\eta z)^{2k}b_q-b=0.
		\end{matrix}
	\end{equation*}
	Thus $b'$ is a theta function of level $2k$, and for $a$ and $d$ to exist we require $[b'x]\in H^1(E_q,\sO)$ be zero, or equivalently that $b'x$ have zero constant term. Since we are working with $\fsl_2$ we also require $a+d=0$. The map $\tau$ is clear, and the calculation of $\delta$ follows the pattern of Lemma \ref{l:extclass}, giving
	\begin{equation*}
		\left[\begin{matrix}
			\tilde{a}_q-\tilde{a} + b'x & 0 \\
			\frac{\tilde{c}_q}{(\eta z)^{2k}}-\tilde{c}- \tilde{a}x+\tilde{d}x+b'x^2 & \tilde{d}_q-\tilde{d}-b'x
		\end{matrix}
		\right]
	\end{equation*}
	where we have chosen various lifts. To clean this up pick $\tilde{c}=0$, and $\tilde{a}=a$ and $\tilde{d}=d=-a$ so that diagonal terms are zero. This gives the simple representative for the cohomology class from the statement.
\end{proof}

\begin{rem}
({\it Skew-symmetry}) We expect the Poisson bivector to be skew symmetric. To check this consider
\begin{equation*}
	\bra B^\#(\omega),\omega\ket=\Tr\left[\left(\begin{matrix}
		0 & 0 \\
		(b'x-2a)x & 0
	\end{matrix}
	\right)\cdot\left(\begin{matrix}
		a & b' \\
		0 & -a
	\end{matrix}
	\right)\right] = b'x(b'x-2a).
\end{equation*}
This is a cohomology representative in $H^1(E_q,\sO)$ and to show it is zero we must check the constant term is zero. Letting $b'x=\sum_l\beta_l z^l$ with $\beta_0=0$ we have $a=\sum \frac{\beta_l}{1-q^l}z^l$ and so
\begin{equation*}
	b'x-2a=\sum_{l\neq0}\left(1-\frac{2}{1-q^l}\right)\beta_lz^l
	= \sum_{l\neq0}\frac{q^l+1}{q^l-1}\beta_l z^l.
\end{equation*}
Then $b'x(2a-b'x)$ has constant term
\begin{equation*}
	\sum_{l\neq0}\frac{q^l+1}{q^l-1}\beta_l\beta_{-l}=
	\sum_{l>0}\left[\frac{q^l+1}{q^l-1} + \frac{q^{-l}+1}{q^{-l}-1}\right]
	\beta_l\beta_{-l} = 0
\end{equation*}
as required.
\end{rem}

\subsection{Calculation of Poisson bracket}
	\label{ss:compare}
	
We use Proposition \ref{p:qbivec} to write out the Poisson bracket on an affine piece of the projective space given by $SL_2$ parabolic bundles lifting a fixed line bundle.

Let $\{[1],\ldots,[z^{2k-1}]\}$ be a basis for $H^1(E_q,\fu(\xi))$. We are interested in the Poisson bracket on the variety $\left(H^1(E_q,\fu(\xi))-\{0\}\right)/\Aut(\xi_0)=\bP^{2k-1}$. Let $x=[x_0+\cdots+x_{2k-1}z^{2k-1}]$ be a point in $H^1(E_q,\fu(\xi))$ with $x_m=x_n=0$.  We calculate the bracket $\bra B^\#_{[x]}(d\vartheta^{2k}_{-m}), d\vartheta^{2k}_{-n}\ket$, where (\ref{e:theta}) 
$\vartheta^{2k}_{-m}=\sum_{l\in\bZ}q^{-ml}q^{kl(l-1)}\eta^{2kl}z^{2kl-m}$ which can be evaluated against $[x]$ as in (\ref{e:evaltheta}). Since we will always be working with $\vartheta$-functions of degree $2k$, we drop the superscript from the notation to reduce clutter. Then
\begin{align*}
	\vartheta_{-m}\cdot x & = \left[\sum_l q^{-ml}q^{kl(l-1)}\eta^{2kl}z^{2kl-m}\right] \left[x_0+\cdots+x_{2k-1}z^{2k-1}\right] \\
	& = \sum_{l\in\bZ,j\in \bZ_{2k}} x_j q^{-ml}q^{kl(l-1)}\eta^{2kl}z^{2kl-m+j} \\
	\mbox{so that }\vartheta_{-m} x-2a & =
	\sum_{l\in\bZ,j\in \bZ_{2k}}x_j\frac{q^{2kl-m+j}+1}{q^{2kl-m+j}-1} q^{-ml}q^{kl(l-1)}\eta^{2kl}z^{2kl-m+j}.
\end{align*}
We will be evaluating against $d\vartheta_{-n}$, so we only consider terms in $(\vartheta_{-m}x-2a)x$ which are of form $z^{2kl+n}$:
\begin{equation*}
	(\vartheta_{-m}x-2a)x = \sum_{l\in\bZ,j,s\in\bZ_{2k}} x_j x_s
	\frac{q^{2kl-m+j}+1}{q^{2kl-m+j}-1}q^{-ml}q^{kl(l-1)}\eta^{2kl}z^{2kl-m+j+s}
\end{equation*}
multiply this by $\vartheta_{-n}$:
\begin{equation*}
	\sum_{j,s,l,t}\frac{q^{2kl-m+j}+1}{q^{2kl-m+j}-1}q^{-ml}q^{-nt}q^{kl(l+1)}q^{kt(t+1)}\eta^{2k(l+t)}z^{2k(l+t)-m+j+s-n}x_jx_s
\end{equation*}
and finally the Poisson bracket is given by looking at the constant term
\begin{equation*}
	\sum_{2k(l+t)+j+s=m+n}\frac{q^{2kl-m+j}+1}{q^{2kl-m+j}-1}q^{-(ml+nt)}q^{k(l(l+1)+t(t+1))}\eta^{2k(l+t)}x_jx_s.
\end{equation*}
Now let $u=2kl-m+j=-(2kt-n+s)$. The above becomes
\begin{equation}
	\label{e:vbrack}
	\bra B^\#_{[x]}(d\vartheta_{-m}), d\vartheta_{-n}\ket = 
	\sum_{u\neq0}\sum_{l,t\in\bZ}\frac{q^u+1}{q^u-1}q^{-(ml+nt)}
	q^{k\left[l(l+1)+t(t+1)\right]} \eta^{2k(l+t)}x_{-(2kl-m-u)}x_{-(2kt-n+u)}.
\end{equation}

\begin{rem}
	Equation (\ref{e:vbrack}) can be written more succinctly as
	\begin{equation}
		\bra B^\#_{[x]}(d\vartheta_{-m}), d\vartheta_{-n}\ket =
		\sum_{u\neq0}\frac{q^u+1}{q^u-1}
		\bra\vartheta_{-m},[xz^{-u}]\ket\cdot
		\bra\vartheta_{-n},[xz^{u}]\ket.
	\end{equation}
\end{rem}

\subsection{Symplectic leaves}

In this section we investigate the symplectic leaves in the moduli space of parabolic bundles. By Corollary \ref{c:leaves} we know these are the fibers of the forgetful map $f:\cM(P)_{\xi_0}\rightarrow \cM(G)$ induced by inclusion $P\hookrightarrow G$.

Focusing on $\SL_2$ we consider the composition $g:H^1(E_q,\fu(\xi))\rightarrow \cM(P)_{\xi_0}\rightarrow\cM(G)$. We adapt the method of Slodowy and Helmke \cite{hs:05} from the case of $\deg\xi_0=2$ to the general case.

\begin{defn}
	Let $V$ be a rank two vector bundle with trivial determinant on $E$. The {\it instability index} of $V$ is the integer
	\begin{equation*}
		i(V)=\max\{\deg\sL\,\vline\,\sL\mbox{ a line bundle with } \Hom(\sL,V)\neq0\}.
	\end{equation*}
\end{defn}

\begin{rem}
	The Atiyah-Bott point of $V$ is $(i(V),-i(V))$.
\end{rem}

The following two lemmas are proven in \cite{hs:05}:

\begin{lem}
	\label{l:degj}
	Let $V$ be a rank two vector bundle on $E$ with trivial determinant and instability index $i(V)=k$. Suppose $L$ is a line bundle of degree $j$, $0<j<k$ and there is a non-trivial map $L\rightarrow V$. Then either
	\begin{equation*}
		V=L\oplus L^*\,\mbox{ or }\, V=\sL\oplus\sL^*\,\mbox{ with }\,
		\deg\sL>j.
	\end{equation*}
\end{lem}

The moduli space of $\SL_2$ bundles with Atiyah-Bott point $(-k,k)$ is isomorphic as a variety to $E$ for $k>0$. If $V$ is semistable and the Atiyah-Bott point is $(0,0)$ then things are a little different.

\begin{rem}
	\label{r:linepairs_unst}
	As a consequence of Lemma \ref{l:degj}, the moduli space of rank two bundles with trivial determinant on $E$ and instability index $k$ can be identified with the collection of pairs $\{(L,L^*)\,|\,\deg L=k\}$ via the map $(L,L^*)\mapsto V:=L\oplus L^*$, and so is isomorphic to $E$.
\end{rem}

\begin{lem}
	\label{l:degzero}
	Suppose $V$ is regular semistable and $L$ is a line bundle of degree zero with a non-trivial map $L\rightarrow V$. Then
	\begin{equation*}
		V\approx\left\{\begin{matrix}
			L\oplus L^* & \mbox{ if } L\neq L^* \\
			L\otimes I_2 & \mbox{ if }L=L^*
		\end{matrix}
		\right.
	\end{equation*}
	where $I_2$ is the Atiyah extension
	$\sO\rightarrow I_2\rightarrow \sO$.
\end{lem}

\begin{rem}
	\label{r:linepairs_st}
	Similar to the unstable case above we have a map from dual pairs of line bundles $(L,L^*)$, now of degree zero, to semistable rank two vector bundles on $E$ with trivial determinant. The map is given by
	\[
	V\mapsto\left\{\begin{matrix}
		L\oplus L^* & \mbox{ if } L\neq L^* \\
		L\otimes I_2 & \mbox{ if } L=L^*
	\end{matrix}
	\right.
	\]
	and is surjective as a consequence of Lemma \ref{l:degzero}. The map is a double cover $E\rightarrow \bP^1$ with four branch points given by the four self-dual line bundles on $E$. In the unstable case (Remark \ref{r:linepairs_unst}) it is an isomorphism $E\rightarrow E$.
\end{rem}

We use the lemmas to describe the closure of the fiber $g^{-1}(V)$.

\begin{lem}
	\label{l:image}
	Let $[x]\in\Ext^1(\xi_0,\xi_0^*)$ be the extension class of $V_{[x]}$, and $L$ be a line bundle with $\deg L = j=i(V_{[x]})$, the instability index of $V_{[x]}$. Then $\Hom(L,V_{[x]})\neq0$ iff $[x]\in\im(ev_a^*)$ for some non-trivial map $a\in\Hom(L,\xi_0)$.
\end{lem}

The map $ev_a^*$ is described during the course of the proof.

\begin{proof}
	Since $\Hom(L,\xi_0^*)=0$ we have exact sequence
	\begin{equation*}
		0\rightarrow\Hom(L,V_{[x]})\rightarrow\Hom(L,\xi_0)
		\xrightarrow{\delta}\Ext^1(L,\xi_0^*).
	\end{equation*}
	It follows that
	\begin{equation*}
		\Hom(L,V_{[x]})\neq0\,
		\mbox{ iff there exists nonzero }a\in\Hom(L,\xi_0)\,
		\mbox{ such that }\delta(a)=0.
	\end{equation*}
	Given $a$, we can construct short exact sequence
	\[
		0\rightarrow L\xrightarrow{a}\xi_0\rightarrow \sO_a\rightarrow0,
	\]
	where $\sO_a:=\cok(a)$.
	Note $\dim\Hom(L,\xi_0)=k-j>0$, where $k=\deg(\xi_0)$. Now consider commutative diagram
	\begin{equation*}
		\xymatrix{
			& & 0\ar[d] & 0\ar[d] \\
			& & \Hom(\xi_0\otimes\sO_a,\xi_0)\ar[d]\ar[r] & 	\Ext^1(\xi_0\otimes\sO_a,\xi_0^*)\ar[d]^{ev_a^*} \\
			& 0\ar[d]\ar[r] & \Hom(\xi_0,\xi_0)\ar[d]^{a^*}\ar[r]^\dd & \Ext^1(\xi_0,\xi_0^*)\ar[d]^{a^*} \\
			0\ar[r] & \Hom(L,V_{[x]})\ar[r] & \Hom(L,\xi_0)\ar[r]^\delta & \Ext^1(L,\xi_0^*)\ar[d] \\
			& & & 0
		}
	\end{equation*}
	Here $\partial(Id_{\xi_0})=[x]$, the extension class of $V_{[x]}$ and $a^*(Id_{\xi_0})=a$. The bottom right rectangle commutes, so $a$ lifts to $\Hom(L,V_{[x]})$ if and only if $a^*([x])=[xa]=0$. The rightmost column is Serre dual to
	\begin{equation*}
		0\rightarrow H^0(L\otimes\xi_0)\xrightarrow{a_*}
		H^0(\xi_0^{\otimes 2})\xrightarrow{ev_a}
		H^0(\xi_0^{\otimes 2}\otimes\sO_a)\rightarrow 0.
	\end{equation*}
	The map $ev_a$ is given by evaluating sections of $\xi_0^{\otimes 2}$ at the points of $E$ in $\sO_a$. Since the column is exact $a$ lifts to a non-zero map in $\Hom(L,V_{[x]})$ iff $[x]\in\im(ev^*_a)$.
\end{proof}

\begin{prop}
	\label{p:fibers}
	Suppose $V$ is an $\SL_2$-bundle on $E$ with $i(V)=j$. Then the closure $\overline{g^{-1}(V)}$ is the cone over a $k-j-1$-dimensional variety embedded in $\Gr^{k-j}\Ext^1(\xi_0,\xi_0^*)$.
\end{prop}

\begin{proof}
	Let $L$ be a line bundle of non-negative degree associated to $V$ as in Remarks \ref{r:linepairs_unst} and \ref{r:linepairs_st}. By Lemmas \ref{l:degj} and \ref{l:degzero} the extension class $[x]$ is in $\overline{g^{-1}(V)}$ iff $\Hom(L,V_{[x]})\neq0$. Lemma \ref{l:image} implies this is true iff $[x]\in\im(ev_a^*)$ for some $0\neq a\in\Hom(L,\xi_0)$.
	
	Note that $\dim\im(ev_a^*)=k-j$  so that $\overline{g^{-1}(V)}$ is a union of $k-j$ dimensional subspaces in $\Ext^1(\xi_0,\xi_0^*)$ varying over morphisms $a$. Each $a$ is determined up to scalar by its zeroes $div(a)$, which form an unordered system of $k-j$ points in $E$, possibly with repetitions. The map $div(a)\mapsto \im(ev_a^*)$ is well-defined and gives a ``generalized linear system''
	\begin{equation*}
		\Sym^{k-j}(E)\rightarrow \Gr^{k-j}\Ext^1(\xi_0,\xi_0^*).
	\end{equation*}
	We are not interested in all possible divisors, but only those which come from maps $a\in\Hom(L,\xi_0)$. Let $\Sigma:\Sym^{k-j}E\rightarrow\Pic^{k-j}E$ be the map induced by addition, and let $X_L:=\Sigma^{-1}(\xi_0\otimes L^*)$. The variety $X_L$ is isomorphic to $\Sym^{k-j-1}E$. Then $\overline{g^{-1}(V)}=\bigcup\{\im(ev_a^*)\,\vline\,div(a)\in X_L\}$: the cone in $\Ext^1(\xi_0,\xi_0^*)$ over $X_L$.
\end{proof}

The dimension of $\overline{g^{-1}(V)}$ is $(k-j-1) + (k-j)$. Quotienting out the $\bC^*$-action, the dimension of a the symplectic leaf is $2(k-j-1)$. Note that finding $\overline{g^{-1}(V)}$ does not automatically give us the leaf, we have to remove the points in the boundary by hand.

\begin{eg}
	Fixing $\deg\xi_0=k$, the maximal unstable case is given by the locus $X=\{[x]\in\Ext^1(\xi_0,\xi_0^*)\,\vline\,i(V_{[x]})=k-1\}$. The line bundle $\underline{\Ext}^1_{E_q}(\xi_0,\xi_0^*)$ has degree $2k$ and determines an embedding $E\hookrightarrow \bP(\Ext^1(\xi_0,\xi_0^*))$. The cone in $\Ext^1(\xi_0,\xi_0^*)$ over this embedding is the elliptic singularity given by taking the total space of degree $-2k$ line bundle $\underline{\Hom}_{E_q}(\xi_0,\xi_0^*)$ and collapsing the zero-section to a point: $[x]=0$.
	
	The symplectic leaves are given by quotienting out the natural $\bC^*$-action on these lines: thus the maximally unstable symplectic leaves form an elliptic curves worth of points, each with trivial Poisson structure.
\end{eg}

\begin{eg}
	The subregular case is when $\deg \xi_0=2$. The unstable locus $X\subset \bA$ corresponds to $j=1$ and is given by collapsing the zero section in the total space of a line bundle on $E$ of degree -4. This is a $\tilde{D}_5$ singularity in the classification. The pre-image of an unstable bundle with Atiyah-Bott point $(-1,1)$ is a line (with origin deleted) in this cone, modding out the $\bC^*$-action gives a point, and these points are the unstable leaves, with trivial Poisson bracket as above.
	
	The stable leaves, where $j=0$, are the first interesting case where we have a symplectic variety of dimension 2. By Proposition \ref{p:fibers} $\overline{g^{-1}(V)}$ is the cone over $E_q$ embedded in $\Gr^2\bC^4$: a rank two vector bundle on $E_q$. Notice however that the Proposition does not give us the leaves; we first need to cut out the unstable locus by hand. 
	
	To do this we need to detect when a map $L\rightarrow V_{[x]}$ from a degree zero line bundle factors through a degree 1 line bundle. Let $L_0$ denote the degree zero bundle and $L_1$ the degree one bundle through which the map factorizes. Consider factorization diagram
	\begin{equation*}
		\xymatrix{
			& & & \xi_0\otimes\sO_{b_1}\ar[dr] \\
			L_0\ar[dr]_{\vartheta_\xi}\ar[rr]^{b_0} & & 
			\xi_0\ar[ur]^{ev_{b_1}}\ar[rr]_{ev_{b_0}} & &
			\xi_0\otimes\sO_{b_0} \\
			& L_1\ar[ur]_{b_1}
		}
	\end{equation*}
	this translates into
	\begin{equation*}
		\xymatrix{
			& \Ext^1(L_1,\xi_0^*)\ar[dr] \\
			\Ext^1(\xi_0\otimes\sO_{b_0}, \xi_0^*)\ar[r]^{ev_{b_0}^*}
			\ar[d] & \Ext^1(\xi_0,\xi_0^*)\ar[r]\ar[u] &
		 	\Ext^1(L_0,\xi_0^*) \\
			\Ext^1(\xi_0\otimes \sO_{b_1},\xi_0^*)\ar[ur]_{ev_{b_1}^*}
		}
	\end{equation*}
and cutting out the unstable locus is reduced to cutting out the images of $ev^*_{b_1}$'s that arise in this way. For each $b_0$ there is a unique such factorization, and the image of $ev^*_{b_1}$ is a line -- the line at infinity. Deleting this and modding out the $\bC^*$-action gives an $\bA^1$-bundle over $E_q$ as the symplectic leaf.
\end{eg}

\section{Loop groups}
	\label{s:loop}
	
The multiplicative quotient description of an elliptic curve $E_q=\bC^*/q^\bZ$ leads to a correspondence between principal $G$-bundles on the curve and $q$-twisted conjugacy classes in the loop group $LG$, as described in \S\ref{ss:qconj}. This provides an alternative method for constructing the moduli space of parabolic bundles, and it becomes reasonable to ask if we can find a loop group interpretation of the Poisson bracket.

Using the theory of $r$-matrices we are able to construct a Poisson structure on $LG$, which passes down to the moduli of parabolic bundles. The construction uses Lu's non-abelian moment map. This is required since the $q$-conjugacy action of the loop group on itself turns out to be a Poisson action rather than a Hamiltonian action. The non-abelian moment map takes values in the dual group of $LG$, which is {\it not} in general equal to $L\fg^*$. Unfortunately, the construction provides not one but many  Poisson structures on the moduli space, depending on the choice of $r$-matrix. In the case of $\SL_2$ we are able to use the explicit formula of \S\ref{s:ranktwo} to find the $r$-matrix which gives rise to the correct Poisson structure, compatible with the bracket defined using sheaf cohomology. In general, there should be a unique $r$-matrix on any $LG$ descending to the Poisson bracket on the moduli space of parabolics defined above.

\subsection{Poisson-Lie groups and $r$-matrices}

Our goal is to provide another construction of the Poisson structure on the moduli space of parabolic bundles. We start by constructing a Poisson structure on the twisted looped group $LG_q$ compatible with $q$-conjugation. The standard way of putting a Poisson structure on a Lie group is via an $r$-matrix. We quickly review the constructions here, for detailed descriptions see the papers \cite{frs:98}, \cite{sts:94}, \cite{sevostyanov:00} and references therein.

\begin{defn}
	An $r$-matrix is an skew-symmetric bilinear operator $r\in\End\fg$
	satisfying the modified classical Yang-Baxter (mCYB) equation
	\begin{equation*}
		\label{e:mcyb}
		[rX,rY]-r([rX,y]+[X,rY])+[X,Y]=0.
	\end{equation*}
\end{defn}

Define maps $r_\pm:\fg\rightarrow \fg$ by $r_\pm=\frac{1}{2}(r\pm Id)$.
The Yang-Baxter implies these maps are homomorphisms. Given $f\in\bC[G]$ let $\nabla_f$, $\nabla'_f\in\fg$ denote the left and right gradients:

\begin{align*}
	\left(\nabla_f(x),\xi\right) & = \left(\frac{d}{dt}\right)_{t=0}f(e^{t\xi}x), \\
	\left(\nabla'_f(x),\xi\right) & = \left(\frac{d}{dt}\right)_{t=0}f(xe^{t\xi})\,\,\mbox{ for }\xi\in\fg.
\end{align*}

Now define Poisson bracket -- the {\it Sklyanin} bracket -- on $G$ by
\begin{equation*}
	\label{e:pbracket}
	\{f,g\}:=\frac{1}{2}(r\nabla'_f,\nabla'_g) -\frac{1}{2}(r\nabla_f,\nabla_g)
\end{equation*}

Decompose $L\fg=L\fn_+\oplus L\fh\oplus L\fn_-$ and let $P_{\bullet}$ denote the respective projections. Introduce $r$-matrix
\begin{equation*}
	\label{e:rmatrix}
	{^\theta}\!r={^{\phi}}\!r=P_{L\fn_+}-P_{L\fn_-}+r^0,\,\, r^0=\frac{1+\phi}{1-\phi}P_{L\fh}.
\end{equation*}
where $\phi:a(z)\mapsto a(qz)$. Under the Sklyanin bracket this gives $LG$ the structure of a Poisson-Lie group. Using Semenov-Tian-Shansky's theory of a twisted Heisenberg double \cite{sts:94} the following bracket on $LG_q$ can be constructed
\begin{equation*}
	\label{e:tbracket}
	\{f,g\}=(r\nabla_f,\nabla_g)+ (r\nabla'_f,\nabla'_g) - 2(r^\phi_+\nabla'_f,\nabla_g) - 2(r_-^\phi\nabla_f,\nabla'_g)
\end{equation*}
where $r_+^\phi=\phi\circ r_+$ and $r_-^\phi=r_-\circ\phi^{-1}$. This gives right Poisson action
\begin{align*}
	\alpha:LG_q\times LG & \rightarrow LG_q \\
	(a,g) & \mapsto (g^\phi)^{-1}\cdot a\cdot g.
\end{align*}

\begin{eg}
	\label{eg:sl2rmatrix}
	For $SL_2$ we can find a compatible Poisson structure explicitly. Our $r$-matrix is a modification of the $r$-matrix defined in \cite{frs:98}; we adapt their construction and ensuing calculations to ensure compatibility with the Poisson bracket on parabolic bundles constructed above. Let $\{E,H,F\}$ be the standard basis in $\fsl_2$ and $\{E_l,H_l,F_l\}$ be the (topological) basis for $L\fsl_2=\fsl_2\otimes \bC((z))$, where for example $H_l=H\otimes z^l$. The $r$-matrix is
	\begin{equation}
		r=\sum_{l\in\bZ}E_l\otimes F_l +\frac{1}{2}\sum_{l\in\bZ}\varphi_l \cdot H_l\otimes H_{-l}, \,\,\,
		\mbox{ with } \varphi_l=\frac{1}{1-q^l}\mbox{ for }l\neq0
		\mbox{ and }\varphi_0=\frac{1}{2}
	\end{equation}
	Note $\varphi_l+\varphi_{-l}=1$. On the tensor product of the two 2-dimensional representations of $\fsl_2((t))$ the $r$-matrix looks like
	\begin{equation*}
		\left(\begin{matrix}
			\varphi\left(\frac{t}{s}\right) & 0 & 0 & 0 \\
			0 & -\varphi\left(\frac{t}{s}\right) & \delta\left(\frac{t}{s}\right) & 0 \\
			0 & 0 & -\varphi\left(\frac{t}{s}\right) & 0 \\
			0 & 0 & 0 & \varphi\left(\frac{t}{s}\right)
		\end{matrix}
		\right),\,\,\,\mbox{ where }\varphi(z)=\frac{1}{2}+\sum_{l\neq0} \frac{z^l}{1-q^l}\mbox{ and }\delta(z)=\sum_{l\in\bZ}z^l.
	\end{equation*}
	The bracket can be calculated as follows. Let
	\begin{equation*}
		L(z)=\left(
		\begin{matrix}
			A(z) & B(z) \\
			C(z) & D(z)
		\end{matrix}
		\right)
	\end{equation*}
	be a functional on $LSL_2$, and define operators $L_1=L\otimes Id$ and $L_2=Id\otimes L$ acting on $\bC^2\otimes \bC^2$ with operator $\sigma$ permuting the two copies of $\bC^2$. The bracket is
	\begin{align*}
		\left\{L_1(z),L_2(w)\right\} & = \frac{1}{2}r^-\left(\frac{w}{z}\right) L_1(z)L_2(z)+ \frac{1}{2}L_1(z)L_2(w)r^-\left(\frac{w}{z}\right) \\
		& - L_1(z)r\left(\frac{wq}{z}\right)L_2(w)+L_2(w)\sigma(r) \left(\frac{zq}{w}\right)L_1(z)
	\end{align*}
	where
	\begin{gather*}
		r^-\left(\frac{w}{z}\right) = r\left(\frac{w}{z}\right)- \sigma(r)\left(\frac{z}{w}\right) =
		\left(\begin{matrix}
			\tau\left(\frac{w}{z}\right) & 0 & 0 & 0 \\
			0 & -\tau\left(\frac{w}{z}\right) & \delta\left(\frac{w}{z}\right) & 0 \\
			0 & -\delta\left(\frac{w}{z}\right) & -\tau\left(\frac{w}{z}\right) & 0 \\
			0 & 0 & 0 & \tau\left(\frac{w}{z}\right)
		\end{matrix}
		\right)
		\\
		\mbox{ and }\tau(z)=\varphi(z)-\varphi\left(\frac{1}{z}\right) = \sum_{l\neq0}\frac{1+q^l}{1-q^l}z^l.
	\end{gather*}
	Using this we find
	\begin{align*}
		\{a_m,a_n\} & = 0 \\
		\{a_m,c_n\} & = \frac{3}{2}a_m c_n- \frac{1}{2}\sum_{l\neq0}a_{m-l}c_{n+l} \\
		\{c_m,c_n\} & = 2\sum_{l\neq0}\frac{1+q^l}{1-q^l}c_{m-l}c_{n+l}
	\end{align*}
	or using shorthand $C(z)=\sum_{\in\bZ}c_lz^l$ and $A(z)=\sum_{l\in\bZ} a_lz^l$
	\begin{align*}
		\{A(z),A(w)\} & = 0 \\
		\{A(z),C(w)\} & = \left(2-\frac{1}{2}\delta\left(\frac{w}{z}\right)\right) A(z)C(w) \\
		\{C(z),C(w)\} & = 2\tau\left(\frac{w}{z}\right) C(z)C(w)
	\end{align*}
	The other brackets can be similarly calculated, but we will not need them in the sequel.
\end{eg}

\subsection{Poisson reduction}

Given the nice relationship between $q$-conjugacy classes and $G$-bundles on $E_q$, it is natural to ask if we can produce parabolic bundles in a similar way. In general the question is rather difficult, simply because of the range of cases which can be considered: recall the moduli space of parabolic bundles is fibered over the associated moduli space of Levi bundles, which can itself be rather complicated. To avoid these difficulties we restrict our attention to the special case where the parabolics are maximal, to the case of Borel subgroups. In this situation the associated Levi is the maximal torus $T$, and so the moduli of Levi bundles is isomorphic to a product of copies of $E_q$.

We perform a Poisson reduction similar to that in \cite{frs:98} and \cite{sevostyanov:00}. There the goal is to perform a deformed Drinfeld-Sokolov reduction to obtain the space of $q$-difference operators, which is an infinite-dimensional Poisson manifold. Our goal is to use similar constructions to produce a Poisson structure on the finite-dimensional space of parabolic bundles on $E_q$.

Following Lu \cite{lu:90} we define a notion of non-abelian moment map for a right Poisson action $\alpha:M\times G\rightarrow M$. Given $X\in\fg$ let $X^r$ denote the corresponding right-invariant 1-form on $G^*$. Also write $\alpha_X$ for the vector field on $M$ generated by the action of $\alpha(e^{tX})$ on $M$.

\begin{defn}
	A map $\mu:M\rightarrow G^*$ is a momentum mapping for the Poisson action $\alpha$ if for all $X\in\fg$
	\begin{equation*}
		\alpha_X=-B_M^\#(\mu^*X^r)
	\end{equation*}
\end{defn}
We will be using Lu's moment map for the $q$-conjugation action $\alpha:LG_q\times LG\rightarrow LG_q$. The following proposition is due to Sevostyanov \cite{sevostyanov:00}
\begin{prop}
	\renewcommand{\labelenumi}{\alph{enumi})}
	\begin{enumerate}
		\item
		Elements $a\in LG_q$ admitting twisted factorization
		\begin{equation*}
			a=a_+^\phi\cdot a_-^{-1}, (a_+,a_-)\in\ LG^*
		\end{equation*}
		form an open dense subset $LG'_q$ in $LG_q$. This factorization is unique in a neighbourhood of the identity and $LG'_q$ is a Poisson submanifold of $LG_q$.
		\item
		The restriction of $\alpha$ to $LG'_q$ has a moment map given by $\mu:a\mapsto (a_+,a_-)$.
		\item
		The induced actions of Poisson-Lie subgroups $LB_\pm$ restricted to $LG'_q$ have moment maps
		\begin{equation*}
			\mu_{LB_\mp}(a)=a_\pm
		\end{equation*}		
	\end{enumerate}
\end{prop}

We apply the proposition and Lu's theory to construct the following diagram of dual pairs:

\begin{equation*}
	\xymatrix{
	& LG_q\ar[dr]^\mu\ar[dl] \\
	LB_-\backslash LG_q & & LB_+ & LT\ar[l]
	}
\end{equation*}

The Levi factor $LT$ is included in $LB_+$ and in fact $LB_+=LT\cdot LN_+$ is a semidirect product. Thus the $q$-twisted dressing action of $LB_-$ on $LB_+$ holds $LT$ invariant. The action of $LB_-$ on $LT$ factors through $LT$: this is $q$-conjugation of $LT$ on itself. The orbits then correspond to $T$-bundles on $E_q$. Let $\cO_{\xi_0}$ be a $LB_-$-orbit inside $T$ corresponding to a $T$-bundle $\xi_0$ on $E_q$. The Poisson reduction $LB_-\backslash\mu^{-1}(\cO_{\xi_0})$ is then the space of $\cM(B)_{\xi_0}$, see \S\ref{ss:casimir}, equipped with a reduced Poisson structure.

Finally, in the example below, the technology developed above is applied to $\SL_2$. We explicitly calculate the reduced Poisson bracket and compare it to the computation in \S\ref{s:ranktwo}.

\begin{eg}
	\label{eg:sl2loop}
	We illustrate the construction in the $SL_2$ case, continuing Examples \ref{eg:sl2flag} and \ref{eg:sl2rmatrix}. Let
	\begin{equation*}
		\cO_{\xi_0}=\left\{\left(\begin{matrix}
			q^l/\eta^k z^k & 0 \\
			0 & \eta^k z^k/q^l
		\end{matrix}
		\right)
		\,\vline\, l\in\bZ
		\right\}
	\end{equation*}
	be the $LT$-orbit in $LB_+$ corresponding to degree $k>0$ line bundle $\xi_0$ on $E_q$. Set
	\begin{equation*}
		\label{e:orblift}
		\tilde{\cU}_{\xi_0}:=\mu^{-1}(\cO_{\xi_0})=\left\{\left(
		\begin{matrix}
			q^l/\eta^k z^k & 0 \\
			x(z) & \eta^k z^k/q^l
		\end{matrix}\right)
		\,\vline\, l\in\bZ\,\mbox{ and }\,x(z)\in\sR={\mathcal Hol}(\bC^*) \right\}.
	\end{equation*}
	Let $\bC[\tilde{\cU}_{\xi_0}]=\bC[\ldots,c_{-1},c_0,c_1,\ldots,a_{-k}]$ be the algebra of functions on $\tilde{\cU}_{\xi_0}$, where $c_i$ picks out the $i^{th}$ coefficient of the power series expansion of $x(z)$. The $LN_-$-action is given by
	\begin{equation*}
		\left(\begin{matrix}
			1 & 0 \\
			z^l & 1
		\end{matrix}
		\right):\left(
		\begin{matrix}
			1/\eta^k z^k & 0 \\
			x(z) & \eta^k z^k
		\end{matrix}
		\right) \mapsto\left(
		\begin{matrix}
			1/\eta^k z^k & 0 \\
			x(z) +\eta^kz^{l+k} - \frac{q^l z^l}{\eta^k z^k} & \eta^k z^k
		\end{matrix} \right).
	\end{equation*}
	Direct computation shows the subalgebra $\bC[\tilde{\cU}_{\xi_0}]^{LN_-}$ of invariant functions is given by
	\begin{gather*}
		\bC[\tilde{\cU}_{\xi_0}]^{LN_-}=\bC[\theta_n^{2k},a_{-k}]
		\,\mbox{ for } n=0,\ldots,2k-1 \\
		\label{e:loopfunctional}
		\mbox{ where } \theta_n^{2k}= \sum_{l\in\bZ}q^{nl}q^{kl^2}a_{-k}^{2kl}c_{2kl+n}.
	\end{gather*}
	To recover the space of parabolic bundles we need to further quotient out the $LT$-action. However we know this results in projective space, so we will have no invariant functions under this action. Instead as in \S\ref{ss:compare} we consider an affine piece. The bracket evaluated at a point $[x]$ is then
	\begin{gather*}
\bra B^\#_{[x]}(d\theta_{-m}), d\theta_{-n}\ket
		= 2\sum_{s\neq0}\sum_{l,j}\frac{1+q^s}{1-q^s} q^{nl+mj} q^{k(l^2+j^2)} \eta^{-2k(l+j)} x_{2kj+m-s} x_{2kl+n+s},
		\\ \mbox{ or in more formal shorthand } \\
		\{\theta_m(z),\theta_n(w)\} = 2\tau\left(\frac{w}{z}\right) \theta_m(z)\theta_n(w)
	\end{gather*}

Comparing with (\ref{e:vbrack}) we see that under the mapping $\vartheta_{-m} \mapsto \theta_{-m}$ the brackets coincide up to a factor of $-2$.
\end{eg}

\begin{rem}
	It should be possible to quantize the Poisson bracket. Semenov-Tian-Shansky \cite{sts:92} defines a twisted Heisenberg double which we believe can be used to quantize the twisted Poisson structure on the loop group. A quantum reduction due to Lu \cite{lu:93} should then provide a quantization of the Poisson bracket on the parabolic moduli space. The obstruction to carrying out this program is finding an $R$-matrix corresponding to the $r$-matrix defined above.
\end{rem}

\section*{Acknowledgements}
I'd like to thank my advisor, V. Ginzburg, for suggesting this area and for many helpful conversations. Partially supported by NSF grant DMS-0401164.

\end{document}